\documentclass[a4paper,12pt]{amsart}
\usepackage{graphics}
\usepackage{amssymb}
\usepackage{amsmath}
\usepackage{latexsym}
\usepackage{amscd}
\newtheorem{Th}{Theorem}
\newtheorem{Prop}[Th]{Proposition}

\newtheorem{Rem}[Th]{\sc Remark}

\newtheorem{Lemma}[Th]{Lemma}
\newcommand{\be}{\begin{eqnarray*}}
\newcommand{\ee}{\end{eqnarray*}}
\newcommand{\pkEF}{{\mathcal P}(^k\!E,F)}

\newcommand{\pkucEF}{{\mathcal P}_{\mathop{\rm uc}\nolimits}(^k\!E,F)}

\newcommand{\pkE}{{\mathcal P}(^k\!E)}

\newcommand{\pk}{{\mathcal P}(^k\!}

\newcommand{\pkwb}{{\mathcal P}_{\mathop{\rm wb}\nolimits}(^k\!}

\newcommand{\pkwbE}{{\mathcal P}_{\mathop{\rm wb}\nolimits}(^k\!E)}
\newcommand{\pkwbEF}{{\mathcal P}_{\mathop{\rm wb}\nolimits}(^k\!E,F)}
\newcommand{\pol}{{\mathcal P}(}

\newcommand{\Proof}{\noindent {\it Proof. }}
\newcommand{\fin}{\hspace*{\fill} $\Box$}
\newcommand{\finesp}{\hspace*{\fill} $\Box$\vspace{.5\baselineskip}}

\newcommand{\N}{\ensuremath{\mathbb{N}}}

\newcommand{\ra}{\rightarrow}
\newcommand{\lra}{\longrightarrow}
\newcommand{\Ra}{\Rightarrow}       
\newcommand{\ptp}{\hat{\otimes}_\pi}
\newcommand{\espv}{\vspace{.5\baselineskip}}

\newcommand{\eps}{\epsilon}

\newcommand{\ufdei}{unconditional finite dimensional expansion of
the identity}

\begin{document}
\title{The polynomial property (V)}

\author[M. Gonz\'alez]{Manuel Gonz\'alez}
\address{Departamento de Matem\'aticas \\
      Facultad de Ciencias\\
      Universidad de Cantabria \\ 39071 Santander (Spain)}
\email{gonzalem@ccaix3.unican.es}
\thanks{The first named author was supported
in part by DGICYT Grant PB 97--0349 (Spain)}

\author[J. M. Guti\'errez]{Joaqu\'\i n M. Guti\'errez}
\address{Departamento de Matem\'atica Aplicada\\
      ETS de Ingenieros Industriales \\
      Universidad Polit\'ecnica de Madrid\\
      C. Jos\'e Guti\'errez Abascal 2 \\
      28006 Madrid (Spain)}
\email{jgutierrez@etsii.upm.es}
\thanks{The second named author was supported
in part by DGICYT Grant PB 96--0607 (Spain)}
\thanks{\hspace*{\fill}\scriptsize file ppv.tex}



\begin{abstract}
Given Banach spaces $E$ and $F$, we denote by $\pkEF$ the space of
all $k$-homogeneous (continuous) polynomials from $E$ into $F$,
and by $\pkwbEF$ the subspace of polynomials which are
weak-to-norm continuous on bounded sets. It is shown that if $E$
has an \ufdei, the following assertions are equivalent: (a) $\pkEF
=\pkwbEF$; (b) $\pkwbEF$ contains no copy of $c_0$; (c) $\pkEF$
contains no copy of $\ell_\infty$; (d) $\pkwbEF$ is complemented
in $\pkEF$. This result was obtained by Kalton for linear
operators. As an application, we show that if $E$ has Pe\l
czy\'nski's property (V) and satisfies $\pkE =\pkwbE$ then, for
all $F$, every unconditionally converging $P\in\pkEF$ is weakly
compact. If $E$ has an \ufdei, then the converse is also true.
\end{abstract}

\maketitle



Given two Banach spaces $E$ and $F$, we denote by $\pkEF$ the
space of all $k$-homogeneous (continuous) polynomials from $E$
into $F$, and by $\pkwbEF$ the subspace of polynomials which are
weak-to-norm continuous on bounded sets. This subspace has been
studied by many authors: see, for instance,
\cite{AHV,AP,FGL,GGwb}. Clearly, every polynomial in $\pkwbEF$
takes bounded sets into relatively compact sets. Observe that
${\mathcal P}(^1\!E,F)= {\mathcal L}(E,F)$, the space of (linear
bounded) operators from $E$ into $F$, and that ${\mathcal
P}_{\mathop{\rm wb}\nolimits}(^1\!E,F)= {\mathcal K}(E,F)$, the
space of compact operators. The spaces ${\mathcal P}(^0\!E,F)$ and
${\mathcal P}_{\mathop{\rm wb}\nolimits}(^0\!E,F)$ may be
identified with $F$.

Kalton studied in \cite{K} the structure of the space ${\mathcal
K}(E,F)$. In the present paper, we obtain versions of his results
for the space $\pkwbEF$, showing that $\pkwbEF$ contains a copy of
$\ell_\infty$ if and only if either $E$ contains a complemented
copy of $\ell_1$ or $F$ contains a copy of $\ell_\infty$. We also
prove that, for $E$ having an \ufdei, the following assertions are
equivalent: (a) $\pkEF =\pkwbEF$; (b) $\pkwbEF$ contains no copy
of $c_0$; (c) $\pkEF$ contains no copy of $\ell_\infty$; (d)
$\pkwbEF$ is complemented in $\pkEF$.

As an application, we prove that, if $E$ has property (V) (the
definitions are given below) and $\pkE =\pkwbE$, then every
$k$-homogeneous unconditionally converging polynomial on $E$ is
weakly compact. If $E$ has an \ufdei, then the converse is also
true.

Throughout, $E$ and $F$ will denote Banach spaces, $B_E$ is the
closed unit ball of $E$ and $S_E$ is the unit sphere of $E$; $E^*$
will be the dual of $E$. The set of natural numbers is denoted by
\N. As usual, $(e_n)$ stands for the unit vector basis of $c_0$.

A formal series $\sum x_n$ in $E$ is {\it weakly unconditionally
Cauchy\/} ({\it w.u.C.,} for short) if, for every $\phi\in E^*$,
we have $\sum|\phi (x_n)|<+\infty$. Equivalent definitions may be
seen in \cite[Theorem~V.6]{Di}. The series is {\it unconditionally
convergent\/} if every subseries converges. Equivalent definitions
may be seen in \cite[Theorem~1.9]{DJT}.

A polynomial $P\in\pkEF$ is {\it unconditionally converging\/}
\cite{FU,FUExt} if, for each w.u.C. series $\sum x_n$ in $E$, the
sequence $\left( P\left(\sum_{i=1}^n x_i\right)\right)_n$ is
convergent in $F$. The space of all unconditionally converging
polynomials is denoted by $\pkucEF$. This class has been very
useful for obtaining polynomial characterizations of Banach space
properties (see \cite{GV}). We say that $P\in\pkEF$ is {\it
(weakly) compact\/} if $P(B_E)$ is relatively (weakly) compact in
$F$. Every weakly compact polynomial is unconditionally
converging. For the general theory of polynomials on Banach
spaces, we refer to
\cite{Din,Mu}.

To each polynomial $P\in\pkEF$ we can associate a unique symmetric
$k$-linear mapping $\hat{P}:E\times\stackrel{(k)}{\ldots}\times
E\to F$ so that $P(x)=\hat{P}(x,\ldots,x)$ and an operator
$T_P:E\ra {\mathcal P}(^{k-1}\!E,F)$ given by
$T_P(x)(y):=\hat{P}(x,y,\stackrel{(k-1)}{\ldots},y)$. It is well
known that $P\in\pkwbEF$ if and only if $T_P$ is compact
\cite[Theorem~2.9]{AHV}.

Denote by ${\mathcal L}_{\mathop{\rm s}\nolimits}(E,{\mathcal
P}(^{k-1}\!E,F))$ the space of all operators $C:E\to {\mathcal
P}(^{k-1}\!E,F)$ such that
$$
(C(x))^\wedge (y_1,\ldots,y_{k-1})=(C(y_1))^\wedge
(x,y_2,\ldots,y_{k-1})\qquad (x,y_1,\ldots,y_{k-1}\in E),
$$
where $(C(x))^\wedge$ stands for the symmetric $(k-1)$-linear
mapping associated to $C(x)$.

\begin{Prop}
The mapping $T:\pkEF\ra {\mathcal L}_{\mathop{\rm
s}\nolimits}(E,{\mathcal P}(^{k-1}\!E,F))$ given by $T(P)=T_P$ is
a surjective linear isomorphism.
\end{Prop}

\Proof
Clearly, $T$ is well defined, linear and injective. Since
$$
\|T_P\|=\sup_{x\in B_E}\|T_P(x)\|=\sup_{x,y\in B_E}\|\hat{P}
(x,y,\stackrel{(k-1)}{\ldots},y)\|\leq\|\hat{P}\|\leq
\frac{k^k}{k!}\,\|P\|
$$
\cite[Theorem~2.2]{Mu}, we have that $T$ is continuous. To see that
it is surjective, take $C\in {\mathcal L}_{\mathop{\rm
s}\nolimits}(E,{\mathcal P}(^{k-1}\!E,F))$, and define $A:E\times
\stackrel{(k)}{\ldots}\times E\to F$ by $A(y_1,\ldots,y_k):=
(C(y_1))^\wedge (y_2,\ldots,y_k)$, and
$P(x)=A(x,\stackrel{(k)}{\ldots},x)$. Then
$$
T_P(x)(y)=A(x,y,\stackrel{(k-1)}{\ldots},y)=C(x)(y)\qquad (x,y\in
E).
$$
Hence, $T_P=C$.\finesp

The subspace of all operators in ${\mathcal L}_{\mathop{\rm
s}\nolimits}(E,{\mathcal P}(^{k-1}\!E,F))$ which are compact
(resp.\ weakly compact) will be denoted by ${\mathcal
K}_{\mathop{\rm s}\nolimits}(E,{\mathcal P}(^{k-1}\!E,F))$ (resp.\
${\mathcal W}_{\mathop{\rm s}\nolimits}(E,{\mathcal
P}(^{k-1}\!E,F))$). Given $C\in {\mathcal K}_{\mathop{\rm
s}\nolimits}(E,{\mathcal P}(^{k-1}\!E,F))$, the symmetry of $C$
easily implies that $C(E)\subseteq {\mathcal P}_{\mathop{\rm
wb}\nolimits}(^{k-1}\!E,F)$.\espv

The proof of \cite[Proposition~5.3]{AS} gives:

\begin{Prop}
\label{compl}
For $n\geq m$, the space ${\mathcal P}(^m\!E,F)$ (resp.\
${\mathcal P}_{\mathop{\rm wb}\nolimits}(^m\!E,F)$) is isomorphic
to a complemented subspace of ${\mathcal P}(^n\!E,F)$ (resp.\
${\mathcal P}_{\mathop{\rm wb}\nolimits}(^n\!E,F)$).
\end{Prop}

\begin{Th}
\label{copylinf}
The space $\pkwbEF$ contains a copy of $\ell_\infty$ if and only
if either $F$ contains a copy of $\ell_\infty$ or $E$ contains a
complemented copy of $\ell_1$.
\end{Th}

\Proof
If $\pkwbEF$ contains a copy of $\ell_\infty$, a fortiori the
space ${\mathcal K}\left( E,{\mathcal P}_{\mathop{\rm
wb}\nolimits}(^{k-1}\!E,F)\right)$ contains it. Therefore
\cite[Theorem~4]{K}, either $E$ contains a complemented copy of
$\ell_1$ or ${\mathcal P}_{\mathop{\rm wb}\nolimits}(^{k-1}\!E,F)$
contains a copy of $\ell_\infty$. Repeating the process, we
conclude that either $E$ contains a complemented copy of $\ell_1$
or ${\mathcal P}_{\mathop{\rm wb}\nolimits}(^0\!E,F)\equiv F$
contains a copy of $\ell_\infty$.

Conversely, if $F$ contains a copy of $\ell_\infty$, since
$F\equiv {\mathcal P}_{\mathop{\rm wb}\nolimits}(^0\!E,F)$ is
isomorphic to a subspace of $\pkwbEF$ (Proposition~\ref{compl}),
we obtain that $\pkwbEF$ contains a copy of $\ell_\infty$. If $E$
contains a complemented copy of $\ell_1$, then $E^*$ contains a
copy of $\ell_\infty$; since $E^*={\mathcal P}_{\mathop{\rm
wb}\nolimits}(^1\!E)$ is isomorphic to a subspace of ${\mathcal
P}_{\mathop{\rm wb}\nolimits}(^1\!E,F)$, which is in turn
isomorphic to a subspace of $\pkwbEF$, we obtain that the latter
contains a copy of $\ell_\infty$.\finesp

The proof of \cite[Lemma~2]{K} yields:

\begin{Lemma}
\label{op2linf}
Assume $E$ is separable, $\pkwbEF$ is complemented in $\pkEF$, and
an operator $\Phi:\ell_\infty\to\pkEF$ is given with the following
properties:

{\rm (a)} $\Phi (e_n)\in\pkwbEF$ for all $n$;

{\rm (b)} the subset $\{ \Phi(\xi)(x):\xi\in\ell_\infty , x\in
E\}\subset F$ is separable.

Then, for every infinite subset $M\subseteq\N$, there exists an
infinite subset $M_0\subseteq M$ with $\Phi(\xi)\in\pkwbEF$ for
all $\xi\in\ell_\infty (M_0)$.
\end{Lemma}

\begin{Lemma}
\label{uncompl}
Suppose $E$ contains a complemented copy of $\ell_1$. Then
$\pkwbEF$ is uncomplemented in $\pkEF$ for every $F$ and $k\in\N$
$(k>1)$.
\end{Lemma}

\Proof
As in \cite[Lemma~3]{K}, we can reduce the problem to the case
$E=\ell_1$.

Fix $v\in S_F$ and define the operator
$$
\Phi:\ell_\infty\lra \pk \ell_1,F)
$$
by
$$
\Phi(\xi)(x)=\sum_{i=1}^\infty\xi_i x_i^k v\quad\mbox{for }\xi=
(\xi_i)_{i=1}^\infty\in\ell_\infty\mbox{ and
}x=(x_i)_{i=1}^\infty\in\ell_1.
$$
Since
$$
\|\Phi(\xi)\|=\sup_{x\in B_{\ell_1}}\left\|\sum_{i=1}^\infty
\xi_i x_i^k v\right\|\leq\|\xi\|\cdot\sup_{x\in B_{\ell_1}}
\sum_{i=1}^\infty \left| x_i^k\right|=\|\xi\| ,
$$
$\Phi$ is continuous (easily, it is even an isometric embedding).

We claim that $\Phi(\xi)\in{\mathcal P}_{\mathop{\rm wb}\nolimits}
(^k\!\ell_1,F)$ if and only if $\xi\in c_0$. Indeed, let
$$
T_{\Phi(\xi)}:\ell_1\lra {\mathcal P}(^{k-1}\!\ell_1,F)
$$
be the associated operator given by
$$
T_{\Phi(\xi)}(x)(y)=\sum_{i=1}^\infty
\xi_i x_i y_i^{k-1}v\qquad\mbox{for $x=(x_i),y=(y_i)\in\ell_1$}.
$$
Since, for $n\neq m$,
\be
\lefteqn{\left\|\left(
T_{\Phi(\xi)}(e_n)-T_{\Phi(\xi)}(e_m)\right)(y)\right\|=}\\
&&\left\| T_{\Phi(\xi)}(e_n-e_m)(y)\right\| =
\left\|\sum_{i=1}^\infty\xi_i(\delta_{in}-\delta_{im})y_i^{k-1}v
\right\|=\left|\xi_ny_n^{k-1}-\xi_my_m^{k-1}\right| ,
\ee
we have
$$
\max\{|\xi_n|,|\xi_m|\}\leq\left\|
T_{\Phi(\xi)}(e_n)-T_{\Phi(\xi)}(e_m)\right\|\leq |\xi_n|+|\xi_m|,
$$
it follows that $T_{\Phi(\xi)}$ is compact if and only if $\xi\in
c_0$ (see \cite[Exercise~VII.5]{Di}).

By Lemma~\ref{op2linf}, there is an infinite $M\subseteq \N$ such
that $\Phi(\xi)\in{\mathcal P}_{\mathop{\rm wb}\nolimits}
(^k\!\ell_1,F)$ for all $\xi\in\ell_\infty (M)$, which contradicts
the above claim.\finesp

In the linear case $(k=1)$, Kalton needs to assume that $F$ is
infinite dimensional for the validity of Lemma~\ref{uncompl}
\cite[Lemma~3]{K}. Taking $k>1$, we can drop this condition.
Observe that, if $k=1$ and $\dim (F)<\infty$, we have
$$
{\mathcal P}_{\mathop{\rm wb}\nolimits}(^1\!E,F)=
{\mathcal K}(E,F)={\mathcal L}(E,F)={\mathcal P}(^1\!E,F),
$$
so the conclusion of the Lemma is not true.

An {\it\ufdei\/} for a Banach space $E$ is a sequence of finite
dimensional operators $A_n:E\to E$ such that for each $x\in E$,
$$
x=\sum_{n=1}^\infty A_n(x)
$$
unconditionally. This condition is a bit more general than having
an unconditional basis \cite{JLS}.

\begin{Lemma}
\label{sumpol}
Suppose $E$ has an \ufdei\ and let $P\in\pkEF$. Then there is a
w.u.C. series $\sum P_i$ in $\pkwbEF$ such that, for all $x\in E$,
$P(x)=\sum_{m=1}^\infty P_m(x)$ unconditionally.
\end{Lemma}

\Proof
There is a sequence $(A_i)\subset {\mathcal K}(E,E)$ such that,
for every $x\in E$, we have $x=\sum_{i=1}^\infty A_i(x)$
unconditionally. Then,
\be
P\left(\sum_{i=1}^n A_i(x)\right)&=&\sum_{i_1,\ldots,i_k=1}^n
\hat{P}\left( A_{i_1}(x),\ldots,A_{i_k}(x)\right)\\
&=&\sum_{m=1}^n\left( \sum_{\max\{ i_1,\ldots,i_k\} =m}\hat{P}
\left(  A_{i_1}(x),\ldots,A_{i_k}(x)\right)\right)\\
&=&\sum_{m=1}^n P_m(x)
\ee
where $P_m\in\pkwbEF$.

Choosing finite subsets $I_1,\ldots,I_k$ of integers, we have
\be
\left\|\sum_{i_1\in I_1,\ldots,i_k\in I_k} \hat{P}\left(
A_{i_1}(x),\ldots,A_{i_k}(x)\right)\right\|&=
&\left\|\hat{P}\left(\sum_{i_1\in I_1} A_{i_1}(x),\ldots,\sum_{i_k\in
I_k}A_{i_k}(x)\right)\right\|\\
&\leq&\|\hat{P}\|\cdot\left\|\sum_{i_1\in I_1}A_{i_1}(x)\right\|
\cdot\ldots\cdot\left\|\sum_{i_k\in I_k}A_{i_k}(x)\right\| .
\ee
Hence, the series
$$
\sum_{i_1,\ldots,i_k=1}^\infty\hat{P}\left(
A_{i_1}(x),\ldots,A_{i_k}(x)\right)
$$
is unconditionally convergent for all $x\in E$
\cite[Theorem~1.9]{DJT}. Therefore, $P(x)=\sum_{m=1}^\infty P_m(x)$
unconditionally.

Moreover, by the uniform boundedness principle
\cite[Theorem~2.6]{Mu}, we have
$$
\sup_{I\subset\N\,\mathop{\rm
finite}}\left\|\sum_{i\in I}P_i\right\|<\infty .
$$
So $\sum P_i$ is w.u.C. in $\pkwbEF$ \cite[Theorem~V.6]{Di}.\fin

\begin{Th}
\label{Kalton}
Suppose $E$ has an \ufdei\ and let $k\in\N$ $(k>1)$. Then the
following assertions are equivalent:

{\rm (a)} $\pkEF =\pkwbEF$;

{\rm (b)} $\pkwbEF$ contains no copy of $c_0$;

{\rm (c)} $\pkEF$ contains no copy of $\ell_\infty$;

{\rm (d)} $\pkwbEF$ is complemented in $\pkEF$.
\end{Th}

\Proof
(a) $\Ra$ (c): Assume $\pkEF =\pkwbEF$ contains a copy of
$\ell_\infty$. By Theorem~\ref{copylinf}, either $E$ contains a
complemented copy of $\ell_1$ or $F$ contains a copy of
$\ell_\infty$. Lemma~\ref{uncompl} implies that $E$ contains no
complemented copy of $\ell_1$, so $F$ contains a copy of
$\ell_\infty$. Take a normalized basic sequence $(x_i)\subset E$
and a bounded sequence of coefficient functionals $(\phi_n)\subset
E^*$ ($\phi_i(x_j)=\delta_{ij}$). Define $P\in\pk E,\ell_\infty)$
by $P(x):=\left(\phi_n(x)^k\right)_n$. Then, for $i\neq j$, we
have $\|P(x_i)-P(x_j)\|\geq |\phi_i(x_i)^k-\phi_i(x_j)^k|=1$.
Hence, denoting by $J:\ell_\infty\ra F$ an isomorphism, we get
$J\circ P\in\pkEF\backslash\pkwbEF$, a contradiction.

(c) $\Ra$ (b): The proof is the same of the linear case
\cite[Theorem~6]{K} with slight modifications.

(b) $\Ra$ (a): Take $P\in\pkEF$. Consider the w.u.C. series $\sum
P_i$ given by Lemma~\ref{sumpol} with $P_i\in\pkwbEF$. Since this
space contains no copy of $c_0$, the series is unconditionally
convergent \cite[Theorem~V.8]{Di}. Clearly, its sum must be $P$.
Since the space $\pkwbEF$ is closed, we conclude that
$P\in\pkwbEF$.

(a) $\Ra$ (d) is trivial.

(d) $\Ra$ (a): If $\pkwbEF$ is complemented in $\pkEF$,
Lemma~\ref{uncompl} implies that $E$ contains no complemented copy
of $\ell_1$. Suppose there is $P\in\pkEF\backslash\pkwbEF$. By
Lemma~\ref{sumpol}, we can find a sequence $(P_i)\subset\pkwbEF$
so that $P(x)=\sum_{i=1}^\infty P_i(x)$ unconditionally for each
$x\in E$, and $\sum P_i$ is w.u.C. but is not unconditionally
convergent since $P\notin\pkwbEF$. Hence, we can find $\eps>0$,
and an increasing sequence $(m_j)$ of integers such that, for each
$j$, the polynomial
$$
C_j:=\sum_{i=m_j+1}^{m_{j+1}}P_i
$$
satisfies $\|C_j\|>\eps$.

Define now $\Phi:\ell_\infty\to\pkEF$ by $\Phi(\xi)(x)=\sum
\xi_jC_j(x)$ for $\xi=(\xi_j)_{j=1}^\infty\in\ell_\infty$ and
$x\in E$. Since the series $\sum P_i(x)$ is unconditionally
convergent, so is the series $\sum\xi_jC_j(x)$. The set
$\{\Phi(\xi)(x):\xi\in\ell_\infty ,x\in E\}$ is contained in the
closed linear span of $\{ C_j(E):j\in\N\}$ which is separable by
the compactness of $C_j$. On the other hand,
$\Phi(e_j)=C_j\in\pkwbEF$ for all $j$. By Lemma~\ref{op2linf},
there is an infinite subset $M\subset\N$ such that
$\Phi(\xi)\in\pkwbEF$ for each $\xi\in\ell_\infty (M)$.

Therefore, for each $\xi\in\ell_\infty (M)$, the series
$\sum\xi_jC_j$ is weak subseries convergent. The Orlicz-Pettis
theorem then implies that it is unconditionally convergent. In
particular,
$$
\lim_{\substack{j\in M\\ j\to\infty}}\|C_j\|=0,
$$
a contradiction.\finesp

Observe that the \ufdei\ is used only in (b) $\Ra$ (a) and (d)
$\Ra$ (a).

In the linear case $(k=1)$, the restriction $\dim (F)=\infty$ is
required for the validity of
Theorem~\ref{Kalton} \cite[Theorem~6]{K}.

\begin{Rem}{\rm
In order to highlight the difference between
Theorem~\ref{copylinf} and Theorem~\ref{Kalton}, let us consider
the spaces $\pk\ell_p,\ell_q)$ and $\pkwb\ell_p,\ell_q)$ for
$1<p,q<\infty$ and $k>1$.

If $kq>p$, then the space $\pk\ell_p,\ell_q)$
contains a copy of $\ell_\infty$. Indeed, define
$J:\ell_\infty\to\pk\ell_p,\ell_q)$ by
$$
J(\xi)(x)=\left(\xi_ix_i^k\right)_{i=1}^\infty .
$$
Since, for all $n$,
$$
|\xi_n|=\| J(\xi)(e_n)\|\leq\| J(\xi)\|=\sup_{x\in B_{\ell_p}} \|
J(\xi)(x)\|\leq\|\xi\| ,
$$
we get that $J$ is an isometric embedding.

However, the space $\pkwb\ell_p,\ell_q)$ contains no copy of
$\ell_\infty$ since it is separable. In fact, it is the norm
closure of the space of all finite type polynomials generated by
the mappings of the form $\phi^n\otimes y$, for
$\phi\in(\ell_p)^*$, $y\in\ell_q$
\cite[Proposition~2.7]{AP}.

On the other hand, if $\xi\in c_0$, then $J(\xi)$ is in the
closure of the finite type polynomials, so
$J(\xi)\in\pkwb\ell_p,\ell_q)$. Hence, this space contains
$J(c_0)$, an isometric copy of $c_0$.

Recall moreover that, if $kq<p$, then the space
$\pk\ell_p,\ell_q)$ is reflexive \cite[4.3]{AF}.}
\end{Rem}

Every polynomial $P\in\pkEF$ has an extension, called {\it the
Aron-Berner extension}, to a polynomial $\overline{P}\in\pk
E^{**},F^{**})$ (see \cite{AB,GGMM,GV}). If $P\in\pkwbE$, then
$\overline{P}\in\pk E^{**})$ is weak-star continuous on bounded
sets \cite{Mo}.

Recall that a Banach space has {\it property (V)}, introduced in
\cite{Pe}, if every unconditionally converging operator on $E$ is
weakly compact. Every $C(K)$ space has property (V) \cite{Pe}.

We shall need the following result:

\begin{Th}
\label{propv}
{\rm \cite{GV}} The following assertions are equivalent:

{\rm (a)} The space $E$ has property (V);

{\rm (b)} for all $F$ and $k\in\N$, the Aron-Berner extension of
every $P\in\pkucEF$ is $F$-valued;

{\rm (c)} There is $k\in\N$ such that, for all $F$, the
Aron-Berner extension of every $P\in\pkucEF$ is $F$-valued.
\end{Th}

Easily, if a polynomial is weakly compact, then its Aron-Berner
extension is $F$-valued \cite{C}.

\begin{Th}
\label{polpropv}
For $k\in\N$, consider the assertions:

{\rm (a)} $E$ has property (V) and $\pkE =\pkwbE$;

{\rm (b)} for each $F$, every $P\in\pkucEF$ is weakly compact.

Then {\rm (a)} $\Ra$ {\rm (b)}. If, moreover, $E$ has an \ufdei,
then {\rm (b)} $\Ra$ {\rm (a)}.
\end{Th}

\Proof
(a) $\Ra$ (b): Given $P\in\pkucEF$, by Theorem~\ref{propv}, the
range of its Aron-Berner extension $\overline{P}$ is contained in
$F$. Take a net $(x_\alpha)\subset B_E$. We can assume that
$(x_\alpha)$ is weak Cauchy and so it converges in the weak-star
topology to some $z\in E^{**}$.

Let $\psi\in F^*$. Then $\psi\circ P\in\pkwbE$ and so,
$$
\psi\circ P(x_\alpha)=\psi\circ\overline{P}(x_\alpha)=
\overline{\psi\circ P}(x_\alpha)\lra\overline{\psi\circ P}(z)=
\psi\circ\overline{P}(z).
$$
Therefore, the net $(P(x_\alpha))$ is weakly convergent to
$\overline{P}(z)\in F$. So, $P(B_E)$ is relatively weakly compact.

(b) $\Ra$ (a): By the comment preceding this Theorem, and by
Theorem~\ref{propv}, (b) implies that $E$ has property (V). Also,
every polynomial in $\pk E,\ell_1)$ is compact. From this, we
obtain that $\pkE$ contains no copy of $c_0$
\cite[Corollary~8]{DL}. A fortiori, $\pkwbE$ contains no copy of
$c_0$. Since $E$ has an \ufdei, we conclude from
Theorem~\ref{Kalton} that $\pkwbE =\pkE$.\finesp

We do not know if the condition on the existence of an \ufdei\ may
be removed from Theorem~\ref{polpropv}. In fact, if (b) is
satisfied, since $\pol ^{k-1}\!E)$ contains no copy of $c_0$ and
$E$ has property (V), we have
\be
&&\pkE ={\mathcal L}_{\mathop{\rm s}\nolimits}\left(E,\pol
^{k-1}\!E)\right) ={\mathcal W}_{\mathop{\rm
s}\nolimits}\left(E,\pol ^{k-1}\!E)\right) ,\\
&&\pkwbE ={\mathcal K}_{\mathop{\rm
s}\nolimits}\left(E,\pol ^{k-1}\!E)\right) .
\ee
So, we only have to show that ${\mathcal K}_{\mathop{\rm
s}\nolimits}\left(E,\pol ^{k-1}\!E)\right) ={\mathcal
W}_{\mathop{\rm s}\nolimits}\left(E,\pol ^{k-1}\!E)\right)$. There
are many conditions on $E$ (apart from the existence of an \ufdei)
that imply this equality \cite{EJ,F},
and it is not known if there are Banach spaces $X$, $Y$ such that
${\mathcal K}(X,Y)\neq {\mathcal W}(X,Y)$ while ${\mathcal
K}(X,Y)$ contains no copy of $c_0$ \cite{E}.

Recall that the condition $\pkE =\pkwbE$ implies that $E$ contains
no copy of $\ell_1$ \cite{G}.

It is proved in \cite{FU} that if $\pkE =\pkwbE$ and $E$ has
property (u) then, for all $F$, every $P\in\pkucEF$ is weakly
compact. Part (a) $\Ra$ (b) of Theorem~\ref{polpropv} is stronger
than the result of \cite{FU}. The latter cannot be applied, for
instance, to the space $c_0\ptp c_0$ which fails property (u)
\cite{L} while it does have property (V) \cite{EH}. For the
definition of property (u) and its relationship to property (V),
see \cite{Pe}. Recall in particular that, if $E$ has property (u)
and contains no copy of $\ell_1$, then $E$ has property (V)
\cite[Proposition~2]{Pe}.


\begin{thebibliography}{99}

\bibitem{AF} R. Alencar and K. Floret, Weak-strong continuity of
multilinear mappings and the Pe\l czy\'nski-Pitt theorem, {\it J.
Math.\ Anal.\ Appl.} {\bf 206} (1997), 532--546.

\bibitem{AB} R. M. Aron and P. D. Berner, A Hahn-Banach extension
theorem for analytic mappings, {\it Bull.\ Soc.\ Math.\ France}
{\bf 106} (1978), 3--24.

\bibitem{AHV} R. M. Aron, C. Herv\'es and M. Valdivia, Weakly
continuous mappings on Banach spaces, {\it J. Funct.\ Anal.} {\bf
52} (1983), 189--204.

\bibitem{AP} R. M. Aron and J. B. Prolla, Polynomial approximation
of differentiable functions on Banach spaces, {\it J. Reine
Angew.\ Math.} {\bf 313} (1980), 195--216.

\bibitem{AS} R. M. Aron and M. Schottenloher, Compact holomorphic
mappings on Banach spaces and the approximation property, {\it J.
Funct.\ Anal.} {\bf 21} (1976), 7--30.

\bibitem{C} D. Carando, Extendible polynomials on Banach spaces.
Preprint.

\bibitem{Di} J. Diestel, {\it Sequences and Series in Banach Spaces,}
Graduate Texts in Math.\ {\bf 92}, Springer, Berlin 1984.

\bibitem{DJT} J. Diestel, H. Jarchow and A. Tonge, {\it Absolutely
Summing Operators,} Cambridge Stud.\ Adv.\ Math. {\bf 43},
Cambridge Univ.\ Press, Cambridge 1995.

\bibitem{Din} S. Dineen, {\it Complex Analysis in Locally Convex
Spaces,} Math.\ Studies {\bf 57}, North-Holland, Amsterdam 1981.

\bibitem{DL} S. Dineen and M. Lindstr\"om, Spaces of homogeneous
polynomials containing $c_0$ or $\ell^\infty$, in: S. Dierolf, S.
Dineen, P. Doma\'nski (eds.), {\it Functional Analysis (Trier,
1994)}, W. de Gruyter, Berlin 1996, 119--127.

\bibitem{E} G. Emmanuele, Answer to a question by M. Feder about
${\mathcal K}(X,Y)$, {\it Rev.\ Mat.\ Univ.\ Complutense Madrid}
{\bf 6} (1993), 263--266.

\bibitem{EH} G. Emmanuele and W. Hensgen, Property (V) of Pe\l
czy\'nski in projective tensor products, {\it Proc.\ Roy.\ Irish
Acad.} {\bf 95}A (1995), 227--231.

\bibitem{EJ} G. Emmanuele and K. John, Uncomplementability of
spaces of compact operators in larger spaces of operators, {\it
Czech.\ Math.\ J.} {\bf 47} (1997), 19--32.

\bibitem{F} M. Feder, On subspaces of spaces with an
unconditional basis and spaces of operators, {\it Illinois J.
Math.} {\bf 24} (1980), 196--205.

\bibitem{FU} M. Fern\'andez-Unzueta, A new approach to
unconditionality for polynomials on Banach spaces. Preprint.

\bibitem{FUExt} M. Fern\'andez-Unzueta, Unconditionally convergent
polynomials in Banach spaces and related properties, {\it Extracta
Math.} {\bf 12} (1997), 305--307.

\bibitem{FGL} J. Ferrera, J. G\'omez and J. G. Llavona, On
completion of spaces of weakly continuous functions, {\it Bull.\
London Math.\ Soc.} {\bf 15} (1983), 260--264.

\bibitem{GGMM} P. Galindo, D. Garc\'\i a, M. Maestre and J. Mujica,
Extension of multilinear mappings on Banach spaces, {\it Studia
Math.} {\bf 108} (1994), 55--76.

\bibitem{GGwb} M. Gonz\'alez and J. M.  Guti\'errez, Factorization
of weakly continuous holomorphic mappings, {\it Studia Math.} {\bf
118} (1996), 117--133.

\bibitem{G} J. M. Guti\'errez, Weakly continuous functions on
Banach spaces not containing $\ell_1$, {\it Proc.\ Amer.\ Math.\
Soc.} {\bf 119} (1993), 147--152.

\bibitem{GV} J. M. Guti\'errez and I. Villanueva, Extensions of
multilinear operators and Banach space properties. Preprint.

\bibitem{JLS} W. B. Johnson, J. Lindenstrauss and G. Schechtman,
On the relation between several notions of unconditional
structure, {\it Israel J. Math.} {\bf 37} (1980), 120--129.

\bibitem{K} N. J. Kalton, Spaces of compact operators, {\it Math.\
Ann.} {\bf 208} (1974), 267--278.


\bibitem{L} F. Lust, Produits tensoriels projectifs d'espaces de
Banach faiblement s\'equentiellement complets, {\it Colloq. Math.}
{\bf 36} (1976), 255--267.

\bibitem{Mo} L. A. Moraes, Extension of holomorphic mappings from
$E$ to $E''$, {\it Proc.\ Amer.\ Math.\ Soc.} {\bf 118} (1993),
455--461.

\bibitem{Mu} J. Mujica, {\it Complex Analysis in Banach Spaces,}
Math.\ Studies {\bf 120}, North-Holland, Amsterdam 1986.

\bibitem{Pe} A. Pe\l czy\'nski, Banach spaces on which every
unconditionally converging operator is weakly compact, {\it Bull.\
Acad.\ Polon.\ Sci.\ S\'er.\ Sci.\ Math.\ Astr.\ Phys.} {\bf 10}
(1962), 641--648.

\end{thebibliography}
\end{document}